\documentclass[11pt]{article}
\usepackage{amsthm,amsmath,amssymb,amsfonts,url,booktabs}
\usepackage{geometry}
\usepackage{hyperref,enumerate}
\usepackage[shortlabels]{enumitem}
\usepackage[babel]{microtype}
\usepackage[english]{babel}
\usepackage[capitalise]{cleveref}
\usepackage{mathabx}
\usepackage{mathtools}

\newtheorem{theorem}{Theorem}[section]
\newtheorem{prop}[theorem]{Proposition}
\newtheorem{lemma}[theorem]{Lemma}
\newtheorem{conj}[theorem]{Conjecture}
\crefname{prop}{Proposition}{Propositions}
\Crefname{prop}{Proposition}{Propositions}
\DeclareMathOperator{\ind}{ind}

\title{Independent sets in regular uniform hypergraphs}
\author{
An Chang\thanks{Center for Discrete Mathematics and Theoretical Computer Science,
Fuzhou University, Fuzhou, China. E-mail: anchang@fzu.edu.cn}
\and
Lusheng Fang\thanks{School of Mathematics and Statistics, Minnan Normal University,
Zhangzhou, China. E-mail: lsfang@mnnu.edu.cn}
\and
Guorong Gao\thanks{Center for Discrete Mathematics and Theoretical Computer Science,
Fuzhou University, Fuzhou, China. E-mail: grgao@fzu.edu.cn}
\and
Caihong Yang\thanks{School of Mathematics and Physics, China University of Geosciences,
Wuhan, China. E-mail: yangch@cug.edu.cn}
}

\begin{document}
\date{}
\maketitle
\begin{abstract}
We prove that, for all integers $r\ge 2$ and $d\ge 1$, every finite simple $d$-regular $r$-uniform hypergraph $H$ on $n$ vertices satisfies
\[
\operatorname{ind}(H)\le \operatorname{ind}(H_d^r)^{n/(rd)},
\]
with equality if and only if $H$ is a disjoint union of copies of $H_d^r$. Here $\operatorname{ind}(G)$ denotes the number of independent sets in $G$, and  $ H_d^r$ is the $d$-regular $r$-partite $r$-graph
on $rd$ vertices obtained by taking all $r$-sets consisting of one vertex
from a fixed $d$-vertex set and one of $d$ pairwise disjoint $(r-1)$-sets
partitioning the remaining vertices. This settles the conjecture of Balogh, Bollob\'as and Narayanan, who proved the bound for quasi-bipartite hypergraphs.
% We prove the conjecture of Balogh, Bollob\'as and Narayanan that a finite simple
% $d$-regular $r$-uniform hypergraph on $n$ vertices has at most
% $[2^{(r-1)d}+(2^d-1)(2^{r-1}-1)^d]^{n/(rd)}$ independent sets,
% and we determine all cases of equality. For $r\geq 3$ and $d\geq 2$, we remove
% the $(r-1)$-sets that lie in $d$ edges and use a conditional moment bound for vertex links to estimate
% the number of independent sets in the remaining hypergraph. The degree of each
% remaining vertex determines how many of the deleted sets form an edge with it.


\end{abstract}

\section{Introduction}\label{intro}

An independent set in a hypergraph is a set of vertices containing no edge.
Write $\ind(H)$ for the number of independent sets of $H$, including the empty
set. Suppose that every edge has $r$ vertices and every vertex lies in $d$
edges. If there are $n$ vertices, then there are $nd/r$ edges.
How should the edges overlap to leave as many vertex subsets independent as possible?

For graphs, the answer comes from the complete bipartite graph $K_{d,d}$.
An independent set of this graph lies entirely in one of its two vertex
classes, so there are $2^{d+1}-1$ choices. Alon \cite{Alon_1991} first
established the asymptotic upper bound $2^{(1/2+o_d(1))n}$ conjectured
by Granville. Kahn \cite{Kahn_2001} proved the sharp bound for regular bipartite
graphs by an entropy argument, and Zhao \cite{Zhao_2010} extended it to all
regular graphs. Thus every $d$-regular graph $G$ on $n$ vertices satisfies
\[
  \ind(G)^{1/n}\leq \ind(K_{d,d})^{1/(2d)}.
\]
Disjoint unions of copies of $K_{d,d}$ attain equality.

For $r\geq 3$, Balogh, Bollob\'as and Narayanan \cite{Balogh_Bollobas_Narayanan_2021} proposed the
following construction. Take pairwise disjoint sets $S_1,\ldots,S_d,T$,
with $|S_i|=r-1$ and $|T|=d$. The edges of $H_d^r$ are all the sets
$S_i\cup\{t\}$ with $1\leq i\leq d$ and $t\in T$. Thus $H_d^r$ has
$rd$ vertices and $d^2$ edges. Balogh, Bollob\'as and Narayanan
\cite{Balogh_Bollobas_Narayanan_2021} give its independent-set count as
\begin{equation}\label{block-count}
 B_{r,d}:=\ind(H_d^r)=2^{(r-1)d}+(2^d-1)(2^{r-1}-1)^d .
\end{equation}
The same construction at $r=2$ gives $K_{d,d}$. They gave the following conjecture.
\begin{conj}\label{main-conj}
For all $r\geq 2$ and $d\in \mathbb{N}$, if $\mathcal{G}$ is a $d$-regular $r$-graph on $n$ vertices, then
$$
\ind(\mathcal{G}) \leq \ind(H^r_d)^{n/rd}.
$$
\end{conj}

Our main result establishes Conjecture
\ref{main-conj} of Balogh, Bollob\'as and Narayanan \cite{Balogh_Bollobas_Narayanan_2021}
and characterises equality.

\begin{theorem}\label{mainthm}
Let $r\geq 2$ and $d\geq 1$ be integers. Every finite simple $d$-regular
$r$-uniform hypergraph $H$ on $n$ vertices satisfies
\[
  \ind(H)\leq B_{r,d}^{\,n/(rd)}.
\]
Equality holds if and only if $H$ is a disjoint union of copies of $H_d^r$.
\end{theorem}

Balogh, Bollob\'as and Narayanan \cite{Balogh_Bollobas_Narayanan_2021}
proved the bound for quasi-bipartite hypergraphs. These admit a vertex partition $A\cup B$
such that each edge meets $A$ once and, for each $a\in A$, the edges through
$a$ are disjoint outside $a$. This includes $H_d^r$ itself, with $A=T$.
For linear hypergraphs without cross-edges, Cohen, Perkins, Sarantis and
Tetali \cite{Cohen_Perkins_Sarantis_Tetali_2022} obtained an upper
bound on $n^{-1}\log\ind(H)$ that is sharp to first order for fixed $r\geq3$
as $d$ tends to infinity. For hypergraphs with controlled
edge overlaps, another approach comes from the
container theorems of Balogh, Morris and Samotij \cite{Balogh_Morris_Samotij_2015} and Saxton
and Thomason \cite{Saxton_Thomason_2015}. Sarantis, Tetali and Zheng
\cite{Sarantis_Tetali_Zheng_2026} used
containers to obtain the conjectured asymptotic exponential rate when, after vertices with identical
incident-edge sets are identified,
the maximum pair codegree is $o(d)$ for fixed $r$. 

For a graph $G$, give each independent set $I$ weight $\lambda^{|I|}$ and write $Z_G(\lambda)$ for the sum of
these weights. The mean proportion of occupied vertices under the resulting
probability law is the \emph{occupancy fraction}. It equals
$\lambda(\log Z_G)'(\lambda)/|V(G)|$, so integrating it recovers the
normalised logarithm of $Z_G$. Galvin and Tetali \cite{Galvin_Tetali_2004} obtained the
weighted bipartite graph bound within a more general result on graph
homomorphisms. Davies, Jenssen, Perkins and Roberts \cite{Davies_Jenssen_Perkins_Roberts_2017} showed
that $K_{d,d}$ maximises the occupancy fraction at every positive activity
among $d$-regular graphs. Their stability argument \cite{Davies_Jenssen_Perkins_Roberts_2018} also gives
a bound for each prescribed size $k>\varepsilon n$. More precisely, for
fixed $d\geq2$ and $\varepsilon>0$, if $n$ is sufficiently large and $2d\mid n$,
the union of $n/(2d)$ copies of $K_{d,d}$ maximises the number of such
independent sets. Davies, Jenssen and Perkins \cite{Davies_Jenssen_Perkins_2021} extended this
conclusion to every $k$ for all sufficiently large multiples $n$ of $2d$
at fixed $d$. For each $r\geq3$ and all sufficiently large degrees $d$, Sarantis, Tetali
and Zheng \cite{Sarantis_Tetali_Zheng_2026} constructed a $d$-regular $r$-uniform hypergraph
whose occupancy fraction at $\lambda=1$ exceeds that of $H_d^r$.
This rules out a direct extension of the pointwise comparison for graphs.

In the case $r\geq 3$ and $d\geq 2$, we first delete every $(r-1)$-set
that lies in $d$ edges. Every edge meeting
such a set contains it, so the choices on its vertices can be counted once
the remaining vertices have been exposed. If a remaining vertex $v$ has
degree $\ell$ after deletion, then exactly $d-\ell$ of the deleted sets form
an edge with $v$. This is why our occupancy estimate keeps track of the
individual degrees in the remaining hypergraph. When every edge contains a
deleted set, the counting bound and its equality condition follow from the
two-activity inequality of Sah, Sawhney, Stoner and Zhao
\cite{Sah_Sawhney_Stoner_Zhao_2019} on the incidence graph of the
deleted sets and their neighbours. When edges remain, the graphical
Brascamp--Lieb inequality of the same authors
\cite{Sah_Sawhney_Stoner_Zhao_2020} allows the vertex weights
in this incidence graph to depend on the remaining degrees.

The \hyperref[mainproof]{proof of Theorem~\ref*{mainthm}}
in Section~\ref{proof} uses two estimates on vertex degrees, stated together
in Proposition~\ref{degree}. Section~\ref{est} proves the proposition.

\section{Preliminary}\label{prelim}
All hypergraphs are finite and simple. For a $r$-graph $G$, write $e(G) =|E(G) |$, $d_{G} (v)$ for the degree of $v$, and
$\Delta_{j} (G)$ for the largest number of edges containing a fixed $j$-set. Put $\Delta (G)
=\Delta_{1} (G)$. The link $L_{G} (v)$ has edges $e \setminus \{v\}$ as $e$ ranges over the edges
containing $v$. For a set $U$ of vertices, $G-U$ denotes the induced hypergraph on $V(G)-U$.
Here $[m] =\{1,\dots ,m\}$ for a positive integer $m$. Throughout the deletion argument, $d$ denotes
the original degree bound. All logarithms are natural.

We first recall the independent-set bound of Kahn and Zhao
\cite{Kahn_2001,Zhao_2010}, including its equality condition;
see also \cite{Sah_Sawhney_Stoner_Zhao_2019}.

\begin{theorem}\label{thm:kahn-zhao}
Let $d\geq 1$ be an integer, and let $G$ be a finite simple
$d$-regular graph on $n$ vertices. Then
\[
\ind(G)\leq \ind(K_{d,d})^{n/(2d)}
          =\bigl(2^{d+1}-1\bigr)^{n/(2d)}.
\]
The bound is attained precisely by disjoint unions of copies of
$K_{d,d}$.
\end{theorem}

For $2$-regular $r$-uniform hypergraphs on $n$ vertices, the bound
$\ind(H)\leq B_{r,2}^{n/(2r)}$ was proved for odd $r\geq 3$ and
for $r=4$ in
\cite{Sarantis_Tetali_Zheng_2026}.
A proof for every $r\geq 2$ is presented in the unrefereed preprint
\cite{Ferudun_2026}. Our proof of \cref{mainthm} treats $d=2$
directly and also determines the equality cases.
The following consequence of the graphical Brascamp--Lieb inequality
\cite{Sah_Sawhney_Stoner_Zhao_2020} will be used
in the reconstruction step.

\begin{lemma}\label{lem:bipartite-bernoulli}
Let $D=(L,R;E)$ be a finite simple bipartite graph, and write
$d_i=d_D(i)$ for $i\in L$ and $m_v=d_D(v)$ for $v\in R$.
Let $(Y_i)_{i\in L}$ and $(Z_v)_{v\in R}$ be jointly independent
Bernoulli random variables with
\[
\Pr(Y_i=0)=q_i,\qquad
\Pr(Z_v=1)=z_v,\qquad
0\leq q_i,z_v\leq 1.
\]
Then
\[
\Pr\bigl(
  Y_i=0\ \text{or}\ Z_v=1
  \text{ for every }iv\in E
\bigr)
\leq
\prod_{iv\in E}
\left[
  q_i^{m_v}+\bigl(1-q_i^{m_v}\bigr)z_v^{d_i}
\right]^{1/(d_i m_v)}.
\]
An empty product is understood to be $1$.
\end{lemma}

\begin{proof}
Apply the graphical Brascamp--Lieb inequality to
$f_{iv}(y,z)=\mathbf{1}_{\{y=0\ \text{or}\ z=1\}}$.
The factor corresponding to $iv$ uses $m_v$ independent copies
of $Y_i$ and $d_i$ independent copies of $Z_v$,
with all these copies jointly independent.
All pairwise conditions hold exactly when all copies of $Y_i$
are zero or all copies of $Z_v$ are one. Their probability is
\[
q_i^{m_v}+\bigl(1-q_i^{m_v}\bigr)z_v^{d_i},
\]
and the graphical Brascamp--Lieb inequality gives the exponent
$1/(d_i m_v)$.
\end{proof}
For $r\geq 3$ and $d\geq 2$, put
\[
\begin{gathered}
 k=r-1,\qquad Q=1-2^{-k},\qquad J=-\log Q,\\
 J_* =\log(4/3),\qquad b=Q^d+(1-Q^d)2^{-d}.
\end{gathered}
\]
Write $\Pi (G) =2^{-|V(G) |} \ind (G)$. This is the probability that a uniformly random
vertex set is independent. By complementation, it is also the probability that a random
set is a vertex cover, that is, that it meets every edge. Isolated vertices do not affect
this probability. Since $2^{r d} b=B_{r,d}$, it suffices to show that $\Pi (H) \leq b^{n / (r d)}$.

A $k$-set of codegree $d$ is called \emph{saturated}. In a $d$-regular hypergraph
its vertices are twins in the terminology of Sarantis, Tetali and Zheng
\cite{Sarantis_Tetali_Zheng_2026}, since they belong to exactly the same edges.
We delete the saturated sets and keep the remaining vertices distinct.
Under product Bernoulli-$p$ measure $\mu_{p}$, where $0<p<1$, each coordinate is selected
independently with probability $p$. Set $q=1-p$. For a nonnegative function $f$ of positive mean,
the tilted measure $\mu_{p,f}$ has density $\frac{f}{\mathbb{E}_{p} f}$. A function is decreasing if
$f(A) \geq f(B)$ whenever $A \subseteq B$. Fixing coordinates or averaging over a subset of
independent coordinates preserves this property.

Write $\lambda =\frac{p}{q}$ and $Z_{G} (\lambda ) =\sum_{I \text{ independent}} \lambda^{|I|}$.
Under the hard-core law an independent set has probability $\frac{\lambda^{|I|}}{Z_{G} (\lambda )}$.
For $N=|V(G)|>0$ write $\alpha_G(\lambda)=\mathbb{E}|I|/N$. The occupancy
identity $\lambda(\log Z_G)'=N\alpha_G$ used by Davies, Jenssen, Perkins and
Roberts \cite{Davies_Jenssen_Perkins_Roberts_2017} gives, after subtracting
the contribution of the unconstrained product measure,
\[
\log \Pi \left(G\right) =N \int_{0}^{1} \left(\alpha_{G} \left(\lambda \right) -\frac{\lambda}{1+\lambda}\right) \frac{\mathrm{d} \lambda}{\lambda}.
\]
All endpoint integrals at zero are understood by continuity. We use $0<p\leq \frac{1}{2}$,
equivalently $0<\lambda \leq 1$, for the integrated estimates.

Harris's correlation inequality \cite{Harris_1960} gives nonnegative covariance
for two increasing functions of independent Bernoulli variables, and likewise
for two decreasing functions. In particular, a decreasing tilt can only raise
the mean of a decreasing function or lower that of an increasing function. A vertex is blocked if at least one edge of its link is fully occupied. We shall
bound the expected number of fully occupied link edges given that the vertex is blocked.
The bound must hold under a decreasing tilt.

Let $e_{1} ,\dots ,e_{m}$ be distinct $k$-sets, where $k,m\geq 1$. For a random set $Y$, put $A_{i}
=\{e_{i} \subseteq Y\}$, $I_{i} =\mathbf{1}_{A_{i}}$ and $C=\sum_{i} I_{i}$. Conditional means given
$C>0$ are used only when that event has positive probability. Write $M_{m} =\frac{m p}{1-q^{m}}$.
This increases with $m$ since $\frac{1+q+\dots +q^{m-1}}{m}$ decreases. For later use, the upper-shadow form of the Kruskal--Katona theorem
\cite{Katona_1968,Kruskal_1963} gives, for every nonempty $S\subseteq[m]$,
\[
  \mu_p\Bigl(\bigcup_{i\in S} A_i\Bigr)\ \geq\ p^{k-1}(1-q^{|S|}).
\]
Cutler and Radcliffe \cite{Cutler_Radcliffe_2013} use the same layerwise
shadow comparison to count independent sets. Here we weight each level by
its $\mu_p$-mass. After adding unused vertices to give a ground set of size at
least $k-1+|S|$, the extremal lexicographic family has a common $(k-1)$-set
and $|S|$ distinct remaining vertices, giving the displayed probability.
We next pass from product measure to a decreasing tilt.

\begin{lemma}\label{r2-4}
Under every decreasing tilt, $\mathbb{E}[C \mid C>0] \leq M_{m}$. More generally, under product
measure let $a=p^{k}$. If $M>0$ and $\mu_{p} (\bigcup_{i \in S} A_{i} ) \geq |S|\frac{a}{M}$ for
every nonempty $S \subseteq [m]$, then $\mathbb{E}_{p,f} C\leq M \mu_{p,f} (C>0)$ for every
decreasing tilt.
\end{lemma}
\begin{proof}
Given the second assertion, the first follows from the shadow bound above, since for $j=|S|$ the ratio $j p^k/[p^{k-1}(1-q^j)]$ is $M_j\leq M_m$.

For the second assertion, put $R=\{i:A_i\text{ occurs}\}$ under product measure.
Our hypothesis is $\Pr_p(R\cap S\ne\varnothing)\geq M^{-1}\mathbb E_p|R\cap S|$
for each $S$, with $M\geq 1$ by the singleton case. Applying
Dughmi's criterion \cite[Theorem~3.6, (c) $\Rightarrow$ (b)]{Dughmi_2020}
to the rank-one uniform matroid on $[m]$ and the induced law of $R$ gives
\[
  \mathbb E_p\max_{i\in R} w_i \ \geq\ M^{-1}\sum_i w_i\mu_p(A_i)
  \qquad (w_i\geq 0),
\]
where the empty maximum is zero.

Take $w_i=\mu_{p,f}(A_i)/a$ and order the indices by nonincreasing weight,
breaking ties arbitrarily. Let $D_i$ be the event that no earlier edge is present.
Conditional on $A_i$ under product measure, the coordinates outside $e_i$ are
still independent. Both $\mathbf{1}_{D_i}$ and $f$ decrease in these coordinates,
so Harris's inequality gives
$\mu_{p,f}(A_i\cap D_i)\geq w_i\mu_p(A_i\cap D_i)$.
The first present edge has maximum weight among those in $R$. Summing over $i$ gives
\[
  \mu_{p,f}(C>0)\geq\mathbb E_p\max_{i\in R}w_i
  \geq M^{-1}\sum_i\mu_{p,f}(A_i)=M^{-1}\mathbb E_{p,f} C .
\]
\end{proof}

\section{Proof of Theorem~\ref*{mainthm}}\label{proof}
After we delete all saturated sets, the maximum $k$-codegree is at most $d-1$.
To bound the probability that a random set covers the remaining edges, we assign
each vertex a weight determined by its degree. For degree weights
$\phi_{0} =0,\phi_{1} ,\dots ,\phi_{d} \geq 0$, put
\[
\Phi \left(G\right) =\sum_{v \in V\left(G\right)} \phi_{d_{G} \left(v\right)} , \quad B=\max_{1\leq j\leq d} \left(\phi_{j} -\phi_{j-1} \right) , \quad \rho_{\ell} =\phi_{\ell} +k \ell B.
\]
The last quantity bounds the change in the potential when a vertex of degree $\ell$ is deleted.

\begin{prop}\label{degree}
For every $k=r-1\geq 2$ and $d\geq 2$ there are degree weights as above with the following properties.
\begin{enumerate}\item Every $r$-uniform $G$ with $\Delta (G) \leq d$ and $\Delta_{k} (G) \leq d-1$
satisfies $\log \Pi (G) \leq -\Phi (G)$.
\item For $0\leq \ell <d$, $\rho_{\ell} <\log 2$. For $0\leq \ell \leq d$,
  \[
\log \left[Q^{d-\ell} +\left(1-Q^{d-\ell} \right) 2^{-d} e^{d \rho_{\ell}} \right] \leq \left(1-\frac{k \ell}{r d}\right) \log b+d \phi_{\ell} ,
\]
  with strict inequality for $\ell >0$.
\end{enumerate}
\end{prop}
We construct the weights in Section~\ref{est}. We first consider the edges that meet
a saturated set. If every edge contains such a set and the hypergraph is regular,
the tail vertices give the quasi-bipartition of \cite{Balogh_Bollobas_Narayanan_2021}.
For maximum degree at most $d$, the following reduction gives both the
bound and the equality condition from the two-activity inequality of Sah,
Sawhney, Stoner and Zhao \cite{Sah_Sawhney_Stoner_Zhao_2019}.

\begin{lemma}\label{r3-2}
In an $r$-uniform hypergraph of maximum degree at most $d$, distinct saturated sets are disjoint.
Every edge meeting one contains it and has its remaining vertex outside all saturated sets. If every
edge contains a saturated set, then $\Pi (H) \leq b^{e(H) / d^{2}}$, with equality precisely for
disjoint unions of copies of $H_{d}^{r}$ and isolated vertices.
\end{lemma}
\begin{proof} Two distinct $k$-sets lie together in at most one edge. A vertex shared by two
saturated sets would therefore have degree at least $2d-1>d$. At a vertex of a saturated set all $d$
incident edges already contain the set. Since $2k>r$, two disjoint saturated sets cannot belong to
one edge.

Write $S_1,\ldots,S_t$ for the saturated sets and
$T_i=\{v:S_i\cup\{v\}\in E(H)\}$. Each $T_i$ has size $d$. Suppose every
edge contains a saturated set. Form the bipartite graph $D$ with left part
$[t]$, right part $\bigcup_i T_i$, and edges $iv$ for $v\in T_i$. Every left
degree is $d$ and each right degree $m_v$ is at most $d$. Each edge of $H$
corresponds to one edge of $D$. If there are no edges, the assertion is immediate.

In a uniformly random vertex set, mark $i$ when every vertex of $S_i$ is selected and
mark $v$ when $v$ is selected. The marks are independent, with probabilities
$1-Q$ and $1/2$, respectively. The set is independent in $H$ precisely when
the marked vertices are independent in $D$. Apply
\cite[Theorem~1.6]{Sah_Sawhney_Stoner_Zhao_2019} with left activity $(1-Q)/Q$
and right activity $1$. After normalising by the total weight of all vertex
subsets, the factor at $iv$ is the $1/(dm_v)$ power of the probability of
independence on a complete bipartite graph with $m_v$ left and $d$ right
vertices. At least one side of this graph must be unmarked. Hence
\[
 \Pi(H)\leq\prod_{iv\in E(D)}h_{m_v}^{1/(dm_v)},
 \qquad h_m=2^{-d}+(1-2^{-d})Q^m .
\]
Let $W$ take value $1$ with probability $2^{-d}$ and $Q$ otherwise.
Then $h_m^{1/m}=\|W\|_m$ is strictly increasing in $m$, since $W$ is
positive and nonconstant. As $h_d=b$, we obtain
$\Pi(H)\leq b^{|E(D)|/d^2}=b^{e(H)/d^2}$, strictly if any $m_v<d$.

At equality $D$ is therefore $d$-regular. The equality condition in the same
theorem forces every component to be complete bipartite, hence a copy of
$K_{d,d}$. Replacing each $i$ by $S_i$ gives the asserted components of $H$,
apart from isolated vertices. The converse follows from \eqref{block-count}
and multiplicativity over components.
\end{proof}
We are now ready to prove Theorem~\ref{mainthm}.
\begin{proof}[Proof of Theorem~\ref{mainthm}]
\phantomsection\label{mainproof}
For $r=2$, the assertion follows from \cref{thm:kahn-zhao},
since $H_d^2=K_{d,d}$. For $d=1$, the edges partition the vertex
set, so $H$ is a disjoint union of copies of $H_1^r$ and
\[
\ind(H)=(2^r-1)^{n/r}.
\]
Suppose henceforth that $r\geq 3$ and $d\geq 2$, and put $k=r-1$.
Choose the degree weights given by \cref{degree}.

Let $S_1,\ldots,S_t$ be all the saturated $k$-sets in $H$.
By \cref{r3-2}, these sets are pairwise disjoint, and every edge
meeting one of them contains it. Put
\[
U=V(H)\setminus\bigcup_{i=1}^t S_i,\qquad F=H[U],
\]
and, for each $i\in[t]$, define
\[
T_i=\{v\in U:S_i\cup\{v\}\in E(H)\}.
\]
If $t=0$, we use $[t]=\varnothing$.
Each $T_i$ has size $d$, and every edge of $H$ outside $F$
is uniquely of the form $S_i\cup\{v\}$ with $v\in T_i$.
Every induced subhypergraph of $F$ has maximum degree at most $d$
and maximum $k$-codegree at most $d-1$: a $k$-set of codegree $d$
in $F$ would have been saturated in $H$ and hence deleted.
Thus \cref{degree} applies to every such subhypergraph.

Put $\ell_v=d_F(v)$ for $v\in U$, and let $M=e(F)$.
Since $H$ is $d$-regular, each $v\in U$ belongs to exactly
$d-\ell_v$ of the sets $T_i$. Counting incidences gives
\[
dt=\sum_{v\in U}(d-\ell_v),\qquad
\sum_{v\in U}\ell_v=rM,\qquad
n=|U|+kt.
\]
For every $U_0\subseteq U$, we have
\[
\Phi(F)-\Phi(F-U_0)
\leq \sum_{v\in U_0}\rho_{\ell_v}.
\]
Indeed, the deleted vertices lose their own weights. Each deleted
edge contains at most $k$ surviving vertices, and there are at most
$\sum_{v\in U_0}\ell_v$ deleted edges.
The loss in weight caused by a unit decrease in degree is at most
$B$. Also, $B\geq 0$, since $\phi_d-\phi_0\geq 0$.

Let $X\subseteq V(H)$ be obtained by selecting each vertex
independently with probability $1/2$. We compute $\Pi(H)$ as
the probability that $X$ is a vertex cover of $H$.
For $A\subseteq[t]$, put
\[
U_A=\bigcup_{i\in A}T_i.
\]
Condition on the collection $A$ of indices $i$ for which
$X\cap S_i=\varnothing$. Since the sets $S_i$ are disjoint,
this collection has probability
$(1-Q)^{|A|}Q^{t-|A|}$. Covering all edges through these sets
forces $U_A\subseteq X$. Once these vertices are selected,
the remaining requirement is that $X\cap(U\setminus U_A)$
cover $F-U_A$. Hence
\[
\Pi(H)
=
\sum_{A\subseteq[t]}
(1-Q)^{|A|}Q^{t-|A|}
\,2^{-|U_A|}\Pi(F-U_A).
\]
Applying the first part of \cref{degree} to $F-U_A$ and using
the preceding deletion bound, we obtain
\[
\Pi(H)
\leq
e^{-\Phi(F)}
\sum_{A\subseteq[t]}
(1-Q)^{|A|}Q^{t-|A|}
\prod_{v\in U_A}
\frac{e^{\rho_{\ell_v}}}{2}.
\]

Set
\[
U_+=\{v\in U:\ell_v<d\},\qquad
m_v=d-\ell_v\quad(v\in U_+).
\]
Let $D$ be the bipartite graph with parts $[t]$ and $U_+$,
where $iv\in E(D)$ if and only if $v\in T_i$.
The degree of $i$ is $d$, and the degree of $v$ is $m_v$.
Take jointly independent Bernoulli random variables
$(Y_i)_{i\in[t]}$ and $(Z_v)_{v\in U_+}$ satisfying
\[
\Pr(Y_i=0)=Q,\qquad
\Pr(Z_v=1)=z_v:=\frac{e^{\rho_{\ell_v}}}{2}.
\]
These are valid probabilities, since \cref{degree} gives
$\rho_{\ell_v}<\log 2$ for $v\in U_+$.
Conditioning on $\{i:Y_i=1\}$ shows that the preceding sum is
\[
\Pr\bigl(
  Y_i=0\ \text{or}\ Z_v=1
  \text{ for every }iv\in E(D)
\bigr).
\]
By \cref{lem:bipartite-bernoulli}, this probability is at most
\[
\begin{aligned}
&\prod_{iv\in E(D)}
\left[
  Q^{m_v}+\bigl(1-Q^{m_v}\bigr)z_v^d
\right]^{1/(dm_v)}
\\
&\hspace{2em}=
\prod_{v\in U_+}
\left[
  Q^{d-\ell_v}
  +\bigl(1-Q^{d-\ell_v}\bigr)
   2^{-d}e^{d\rho_{\ell_v}}
\right]^{1/d}.
\end{aligned}
\]
For $\ell_v=d$, the expression in square brackets is $1$.
We may therefore extend the last product to all $v\in U$.
Taking logarithms gives
\[
\log\Pi(H)
\leq
\sum_{v\in U}
\left[
-\phi_{\ell_v}
+\frac1d\log\left(
Q^{d-\ell_v}
+\bigl(1-Q^{d-\ell_v}\bigr)
2^{-d}e^{d\rho_{\ell_v}}
\right)
\right].
\]
By the second part of \cref{degree},
\[
\log\Pi(H)
\leq
\frac{\log b}{d}
\sum_{v\in U}
\left(1-\frac{k\ell_v}{rd}\right).
\]
The incidence identities above imply
\[
\sum_{v\in U}
\left(1-\frac{k\ell_v}{rd}\right)
=
|U|-\frac{k}{rd}\sum_{v\in U}\ell_v
=
\frac{|U|+kt}{r}
=
\frac nr.
\]
It follows that
\[
\Pi(H)\leq b^{n/(rd)}.
\]

If $F$ has an edge, then $\ell_v>0$ for some $v\in U$,
and the reconstruction inequality in \cref{degree} is strict
for that vertex. Thus the final bound is strict in this case.
If $F$ has no edges, every edge of $H$ contains a saturated set,
and \cref{r3-2} gives the equality cases. Since $H$ is
$d$-regular with $d\geq 2$, it has no isolated vertices,
so equality holds precisely when $H$ is a disjoint union
of copies of $H_d^r$.

Finally, $2^{rd}b=B_{r,d}$ gives
\[
\ind(H)=2^n\Pi(H)\leq B_{r,d}^{n/(rd)}.
\]
Conversely, independent-set counts multiply over components,
so every disjoint union of copies of $H_d^r$ attains this bound.
\end{proof}
\section{Estimates}\label{est}
To prove \cref{degree}, we integrate a bound on the occupancy fraction that retains
a separate term for each vertex of positive degree below $d$. For small degrees these extra terms allow us to use quadratic weights in the
reconstruction inequality.

\subsection{Link moments}\label{local}
The Kruskal--Katona bound in Section~\ref{prelim} is attained by a family
with a common $(k-1)$-set. When no such set exists, we obtain a stronger
estimate by applying the same bound separately to two subfamilies.

\begin{lemma}\label{r4-1}
For $d\geq 3$ distinct $k$-sets with no common $(k-1)$-set,
\[
\mu_{p} \left(C>0\right) \geq p^{k-1} \left(1-q^{d-1} \right) +q p^{k} .
\]
The bound is sharp for every $0<p<1$.
\end{lemma}
\begin{proof} Remove the intersection of all the edges, of size $a\leq k-2$.
It is enough to prove the
claim for an $s$-uniform family with empty intersection, where $s=k-a\geq 2$, and multiply by
$p^{a}$. If a vertex has degree $D\geq 2$, set $h=d-D$. Condition on this vertex. When absent,
consider the $h$ avoiding edges. When present, consider the $D$ distinct $(s-1)$-sets in its link.
Applying the shadow bound of Section~\ref{prelim} to each family gives
\[
\begin{aligned}
 \mu_p(C=0)&\leq q[1-p^{s-1}(1-q^h)]+p[1-p^{s-2}(1-q^D)]\\
 &=1-p^{s-1}(1+q-q^{h+1}-q^D).
\end{aligned}
\]
The required bound now follows from $2\leq D\leq d-1$ and
\[
 q^2+q^{d-1}-q^D-q^{d+1-D}
  =(1-q^{D-2})(q^2-q^{d+1-D})\geq0 .
\] If the family is a matching, identify a vertex of one edge with a vertex of another. The resulting edges are distinct because $s\geq2$. The probability of avoiding the two increases by $q p^{2s-1}$, and avoidance of the remaining edges is independent. The modified family satisfies the preceding case. For equality, take a $(d-1)$-edge star in a graph and join two leaves. Adjoin a fixed set of $k-2$ vertices to each edge. Conditioning on the centre gives the claimed value. \end{proof}

For $k\geq 2$ and $0<p\leq \frac{1}{2}$, define the polynomials
\[
\alpha_{k} =\sum_{s=1}^{k-1} s p^{s} , \quad \omega_{k} =k+\sum_{s=1}^{k-1} \left(k-s\right) p^{s} , \quad b_{k} =\sum_{s=1}^{k-1} \frac{s\left(k-s\right)}{k-1}p^{s} .
\]
For a link $L$ let $C_{j} =\sum_{e \in L} \binom{|e \cap Y|}{j}$ so that $C=C_{k}$. Use the score
\[
K=\sum_{j=1}^{k} \frac{k}{\binom{k}{j}} p^{k-j} C_{j} .
\]
A single edge with $k-2,k-1,k$ selected vertices contributes $\alpha_{k} -b_{k} ,\alpha_{k}
,\alpha_{k} +\omega_{k}$, respectively. Its score is increasing with the number selected. Thus every
$\ell$-edge link has conditional mean score at most $\mathcal{U}_{k,\ell} =\ell \alpha_{k}
+\omega_{k} M_{\ell}$.

Conditioning on vertices outside a link can give some of its edges probability zero.
This is why the second part of Lemma~\ref{r2-4} requires a lower bound on the union
probability for every subfamily. The preceding estimate gives these bounds, so we can
apply the lemma under the conditioning in Section~\ref{occ}.

\begin{lemma}\label{r4-2}
If $d\geq 2$ distinct $k$-sets have no common $(k-1)$-set, every decreasing tilt satisfies
$\mathbb{E}[K \mid C>0] \leq \mathcal{K}_{d}$, where
\[
\mathcal{K}_{d} =\begin{cases}2 \alpha_{k} +\frac{\omega_{k} -b_{k} q^{2}}{1-\frac{p^{2}}{2}} \quad & d=2 \\ d \alpha_{k} +\omega_{k} \widehat{M}_{d} \left(p\right) \quad & d\geq 3\end{cases},
\]
\[
\widehat{M}_{3} \left(p\right) =\frac{3}{3-2p}, \quad \widehat{M}_{4} \left(p\right) =\frac{4}{\left(2-p\right)^{2}}, \quad \widehat{M}_{d} \left(p\right) =M_{d-1} \quad \left(d\geq 5\right) .
\]
The conditional mean of $C$ is at most $\widehat{M}_{d}$, with $\widehat{M}_{2} =\frac{2}{2-p^{2}}$.
Each mean-count bound is sharp for every $k$ and $p$ in the stated range. Moreover
$\mathcal{U}_{k,\ell} \leq \mathcal{K}_{d}$ for $1\leq \ell <d$.
\end{lemma}
\begin{proof} \cref{r2-4,r4-1} bound the conditional mean of $C$ by the largest ratio $\frac{m
p^{k}}{\mu_{p} (\bigcup_{i \in S} A_{i} )}$ over nonempty index sets $S$ of size $m$. For a proper
subfamily this is at most $M_{d-1}$. For the full family with $d\geq 3$ it is at most $\frac{d
p}{1-q^{d-1} +p q}$. For $d=3,4$ the latter equals the displayed $\widehat{M}_{d}$ and exceeds
$M_{d-1}$ by direct subtraction. For $d\geq 5$ it is at most $M_{d-1}$, since $q M_{d-1} \geq q
M_{4} \geq 1$ on $q\geq \frac{1}{2}$. The last inequality follows from $4q-(1+q+q^{2} +q^{3} )
=(1-q) (q^{2} +2q-1) \geq 0$. Combining this with $K\leq d \alpha_{k} +\omega_{k} C$ proves the
score bound for $d\geq 3$ and gives the comparison for lower degrees.

For $d=2$, inclusion-exclusion gives $\mu_{p} (C>0) \geq 2p^{k} -p^{k+2}$. Each one-edge subfamily
has ratio one, so \cref{r2-4} gives $\mathbb{E}[C \mid C>0] \leq \frac{2}{2-p^{2}}$. For sharpness
at $d=2$, take two edges intersecting in exactly $k-2$ vertices under product measure. For $d=3,4$
use the example in \cref{r4-1}. For $d\geq 5$ take $d-1$ edges with a common $(k-1)$-set and an
extra edge disjoint from them. Force one vertex of the extra edge absent. Conditional on $C>0$ the
remaining count has law $\operatorname{Bin}(d-1,p)$ conditioned to be positive.

To refine the two-edge score, let $Z$ count edges missing at least two vertices. Write $\omega
=\omega_{k}$, $D=b_{k}$ and $a=\frac{\omega -D q^{2}}{1-\frac{p^{2}}{2}}>0$. Positivity follows by
comparing the coefficients of $D$ and $\omega$. Fix one edge present. The other has $s\geq 2$
unfixed vertices. Under the remaining product measure,
\[
\mathbb{E}\left[\frac{1}{C}\right] =1-\frac{p^{s}}{2}\geq 1-\frac{p^{2}}{2}, \quad \mathbb{E}\left[\frac{Z}{C}\right] =\Pr \left(\operatorname{Bin}\left(s,q\right) \geq 2\right) \geq q^{2} .
\]
The two bounds give $\mathbb{E}[\frac{a+D Z}{C}] \geq a(1-\frac{p^{2}}{2}) +D q^{2} =\omega$. With
the edge fixed, the ratio is a decreasing function of the remaining coordinates, so Harris preserves
this lower bound under a decreasing tilt. Multiply by the probability that the fixed edge is present
and sum over both edges. This gives $\mathbb{E}[(a+D Z) \mathbf{1}_{C>0} ] \geq \omega \mathbb{E}
C$. Since $K\leq 2 \alpha_{k} +\omega C-D Z$, the conditional score is at most $2 \alpha_{k} +a$.
Finally $b_{k} \leq \alpha_{k}$ and $b_{k} <\omega_{k}$, so $2 \alpha_{k} +a>2 \alpha_{k}
+\omega_{k} -b_{k} q^{2} \geq \alpha_{k} +\omega_{k} =\mathcal{U}_{k,1}$. \end{proof}

\subsection{Occupancy}\label{occ}
We use the occupancy method of Davies, Jenssen, Perkins and Roberts
\cite{Davies_Jenssen_Perkins_Roberts_2017}. For hypergraphs without cross-edges,
Cohen, Perkins, Sarantis and Tetali \cite[Sections~2 and~3]{Cohen_Perkins_Sarantis_Tetali_2022}
compare occupancy at a random root and in its neighbourhood in the linear case.
Sarantis, Tetali and Zheng \cite[Section~3.2]{Sarantis_Tetali_Zheng_2026}
make a related comparison after identifying twins. Here the remaining hypergraph
can be irregular and link edges can overlap. We count occupied subsets of each
edge according to their size when averaging over the root in \cref{r4-4}. Fix an
$r$-uniform $G$ with $\Delta(G)\leq d$ and $\Delta_k(G)\leq d-1$. Under the hard-core
law of Section~\ref{prelim}, we first bound from below the probability that every
vertex of a given link edge is occupied.

\begin{lemma}\label{r4-3}
Suppose $\Delta (G) \leq d$ and $\Delta_{k} (G) \leq d-1$. For every $k$-set $S$ contained in an edge,
\[
\Pr \left(S \subseteq I\right) \geq \Lambda_{r,d} \left(p\right) :=\frac{p^{k} T}{1-p^{k} +p^{k} T}, \quad T=q^{d-1} \left(1-p^{2} \right)^{2} .
\]
\end{lemma}
\begin{proof} Let $a_{j}$ count edges meeting $S$ in $j$ vertices. Write $a_{k} =d-m$ with $m\geq
1$. The sum of degrees in $S$ gives $\sum_{j<k} j a_{j} \leq k m$. Draw a hard-core independent set
of $G-S$ and let $t$ be the probability that all of $S$ can be added. Under product measure both the
original independence event and the extension constraints are decreasing. Harris gives
\[
t\geq q^{d-m} \prod_{j=1}^{k-1} \left(1-p^{k+1-j} \right)^{a_{j}} .
\]
For $1\leq j<k$, $p^{k-1-j} \leq 2^{-(k-1-j)} \leq \frac{j}{k-1}$, since $2^{u} \geq 1+u\geq
\frac{k-1}{j}$ at $u=k-1-j$. Concavity at one gives $1-p^{k+1-j} \geq (1-p^{2} )^{\frac{j}{k-1}}$.
Thus
\[
t\geq q^{d-m} \left(1-p^{2} \right)^{k \frac{m}{k-1}} \geq q^{d-m} \left(1-p^{2} \right)^{2m} \geq T,
\]
using $(1-p^{2} )^{2} -q=p(1-2p+p^{3} ) \geq 0$. Decompose $Z_{G}$ by its intersection with $S$. The
extension probability is $t$ for the whole set and at most one for each proper subset. Dividing by
$(1+\lambda )^{k} Z_{G-S}$ in the probability ratio gives $\Pr (S \subseteq I) \geq \frac{p^{k}
t}{1-p^{k} +p^{k} t}\geq \Lambda_{r,d}$. \end{proof}

Set $\mathcal{K}_{\ell} =\mathcal{U}_{k,\ell}$ for $1\leq \ell <d$ and use $\mathcal{K}_{d}$ from \cref{r4-2}. Define
\[
x=\frac{k d p^{k}}{\mathcal{K}_{d} +k d p^{k}}, \quad W_{\ell} =1-\frac{\mathcal{K}_{\ell} +k \ell p^{k}}{\mathcal{K}_{d} +k d p^{k}},
\]
\[
c_{d}^{\ast} = \frac{r}{d} \int_{0}^{1/2} \frac{x}{q} \,\mathrm{d} p, \quad \Delta_{d,\ell} =\int_{0}^{1/2} \frac{W_{\ell} \ell \Lambda_{r,d}}{M_{\ell} q} \,\mathrm{d} p, \quad \Delta_{d,d} =0.
\]
In particular $W_{\ell} \geq 0$ and $\Delta_{d,\ell} \geq 0$.
At a root of degree $d$, Lemma~\ref{r4-2} bounds the mean score given that
the root is blocked. At a root of smaller positive degree, Lemma~\ref{r4-3}
bounds the probability that a given link edge is fully occupied. After averaging
over the root and integrating, we obtain a separate term for each such vertex.

\begin{lemma}\label{r4-4}
Every $G$ satisfying the hypotheses of \cref{r4-3} has
\[
\log \Pi \left(G\right) \leq -c_{d}^{\ast} e\left(G\right) -\sum_{v:1\leq d_{G} \left(v\right) <d} \Delta_{d,d_{G} \left(v\right)} .
\]
\end{lemma}
\begin{proof} Fix a root $v$ and condition on the independent set outside its closed neighbourhood.
Let $Y$ be the selected neighbours and $z=\mathbf{1}_{v \in I}$. With $v$ absent, the allowed $Y$
form a downset, including restrictions from edges crossing the boundary. Given $Y$, the mean of $z$
is $p \mathbf{1}_{C=0}$. Given $C>0$, the root is absent. Under this conditioning, the law of $Y$ is
the Bernoulli-$p$ product law tilted by the downset indicator and conditioned on $C>0$. This is the
law covered by \cref{r4-2}.

Put $a=\mathbb{E}[z]$, $b_{j} =\frac{\mathbb{E}[C_{j} ]}{d \binom{k}{j}}$ and $u_{j}
=\frac{\mathbb{E}[z C_{j} ]}{d \binom{k}{j}}$. Consider
\[
S_{v} =a+x\left(b_{1} -a\right) +\sum_{j=2}^{k} \frac{x}{p^{j-1}}\left(b_{j} -u_{j-1} \right) -p\left(1-x\right) .
\]
For $j<k$, conditioning on $Y$ gives
\[
\frac{b_{j}}{p^{j-1}}-\frac{u_{j}}{p^{j}}=\frac{\mathbb{E}\left[C_{j} \mathbf{1}_{C>0} \right]}{d \binom{k}{j}p^{j-1}}.
\]
Also $a-p=-p \Pr (C>0)$ and $C_{k} =C$ vanishes on $C=0$. Substituting the score and using $p(1-x)
=\frac{x \mathcal{K}_{d}}{k d p^{k-1}}$ gives
\[
S_{v} =\frac{x}{k d p^{k-1}} \mathbb{E}\left[\left(K-\mathcal{K}_{d} \right) \mathbf{1}_{C>0} \right] .
\]
At degree $d$ this is nonpositive: any common $(k-1)$-set in the link, together with $v$, would give
$k$-codegree $d$. At positive degree $\ell <d$ it is at most $-\eta_{\ell} \Pr (C>0)$, where
$\eta_{\ell} =\frac{x(\mathcal{K}_{d} -\mathcal{K}_{\ell} )}{k d p^{k-1}}$. At degree zero it
vanishes.

Now average over the outside configuration and a uniform root. Each occupied $(j+1)$-set in an
original edge is counted $j+1$ times in $z C_{j}$ and $k-j$ times in $C_{j+1}$. Since
$\frac{j+1}{\binom{k}{j}}=\frac{k-j}{\binom{k}{j+1}}$, the averages of $u_{j}$ and $b_{j+1}$ agree.
Write $a_{v} =\Pr (v \in I)$ and $\ell_{v} =d_{G} (v)$. Counting a single occupied vertex gives
$\mathbb{E}[b_{1} ] =\frac{\sum_{v} \ell_{v} a_{v}}{d N}$. Thus, with $h_{v} =\Pr (C_{v} >0)$ and
$p-a_{v} =p h_{v}$, the averaged estimate is
\[
\left(1-x\right) \left(\alpha_{G} -p\right) +\frac{x}{d N}\sum_{v} \ell_{v} a_{v} \leq -\frac{1}{N} \sum_{v:1\leq \ell_{v} <d} \eta_{\ell_{v}} h_{v} .
\]
Using $\sum \ell_{v} =r e(G)$ and $W_{\ell} =x(1-\frac{\ell}{d}) +\frac{\eta_{\ell}}{p}$, we obtain
\[
\alpha_{G} \leq p-\frac{r e\left(G\right)}{d N}x p-\frac{1}{N} \sum_{v:1\leq \ell_{v} <d} W_{\ell_{v}} \left(p-a_{v} \right) .
\]
At degree $\ell$ \cref{r2-4} gives $\mathbb{E}[C_{v} ] \leq M_{\ell} h_{v}$ after averaging over
local views. In the full law \cref{r4-3} gives $\mathbb{E}[C_{v} ] \geq \ell \Lambda_{r,d}$.
Consequently $p-a_{v} \geq \frac{p \ell \Lambda_{r,d}}{M_{\ell}}$. Integrate using $\frac{\mathrm{d}
\lambda}{\lambda}=\frac{\mathrm{d} p}{p q}$. If $N=0$ the assertion is immediate. \end{proof}

It remains to choose degree weights satisfying \cref{degree}. We use
\[
\phi_{\ell} =\frac{c \ell}{r}+\kappa \ell \left(d-\ell \right) , \quad B=\frac{c}{r}+\kappa \left(d-1\right) , \quad \rho_{\ell} =c \ell +\kappa \ell \left(r d-k-\ell \right) ,
\]
with
\[
\left(c,\kappa \right) =\begin{cases}\left(-\frac{\log b}{d^{2}}+\frac{r J^{2}}{90J_{\ast} d^{2}}, \frac{\left(4d+1\right) J^{2}}{60J_{\ast} d 2^{d}}\right) \quad & 2\leq d\leq 9 \\ \left(\frac{J}{d} \left(1+\frac{1}{10d}\right) ,0\right) \quad & d\geq 10\end{cases}.
\]
Indeed $B$ is the first difference at $\ell =1$, the largest by concavity. By \cref{r4-4} and $r
e(G) =\sum_{v} d_{G} (v)$, it suffices for the first part of \cref{degree} that
\begin{equation}\label{target}
\frac{c_{d}^{\ast} \ell}{r}+\Delta_{d,\ell} >\phi_{\ell} \quad \left(1\leq \ell \leq d\right) .
\end{equation}
To check \eqref{target} for all edge sizes, we compare the integral at different ranks.

\subsection{Rank comparison}\label{rank}
Put $J_{k} =-\log (1-2^{-k} )$. The change of variable $t=-\frac{\log (1-p^{k} )}{J_{k}}$ takes
$[0,\frac{1}{2}]$ to $[0,1]$. It allows us to compare the normalised rate
$\frac{c_{d}^{\ast}}{J_{k}}$ at different ranks. The next lemma uses real $k\geq 2$.

\begin{lemma}\label{r4-5}
Fix $0<t\leq 1$. Put $p=p_{k} (t) =(1-e^{-t J_{k}} )^{\frac{1}{k}}$ and
\[
h=h_{k} \left(p\right) =\frac{1}{k+1}\left(\frac{k}{1-p^{k}}-\frac{p}{1-p}\right) .
\]
As $k$ increases, $p$ is nondecreasing and $h$ is strictly increasing. For $t<1$ the product
$g=h(\frac{1}{p}-2)$ is strictly decreasing. At $t=1$ it vanishes.
\end{lemma}
\begin{proof} Put $L=\log 2$, $a=2^{-k}$, $J=J_{k}$, $r=k+1$, $z=p^{k}$ and $u=\frac{p}{q}$. We use
$\frac{2}{3}<L<\frac{7}{10}$ from the series for $2 \operatorname{arctanh}(\frac{1}{3})$ and
$e^{\frac{7}{10}}$. The expansion $1-(1-x)^{t} =\sum_{n\geq 1} c_{n} x^{n}$ has $c_{n} \geq 0$ and
$\sum c_{n} =1$. Hence $p$ is the $L^{k}$ norm of a variable taking the values $2^{-n}$ with
probabilities $c_{n}$. Its monotonicity follows from the power-mean inequality.

Let primes denote total derivatives at fixed $t$, and set $\theta =\frac{t a}{(1-a) (e^{t J} -1)}$.
Convexity of the exponential gives $\theta \geq 1$. Differentiation of $p^{k} =1-e^{-t J}$ gives
$\frac{p'}{p}=\frac{-\log p-L \theta}{k}\leq \frac{\log (\frac{1}{2p})}{k}$. By concavity in $t$,
$z\geq t a$. Since $h=\frac{k e^{t J} -u}{r}$, we have $r^{2} h' =e^{t J} +u+k r e^{t J} t J' -r
u'$. At fixed $p\leq \frac{1}{2}$, the function $k r \frac{p^{k}}{(1-2^{-k} ) (1-p^{k} )}$ decreases
with $k$. Its logarithmic derivative is at most $\frac{1}{k}+\frac{1}{r}-\frac{L(1+a)}{1-a}<0$: use
the increase of $\frac{\log p}{1-p^{k}}$ in $p$ and $x \coth (x) >1$ at $x=k \frac{L}{2}$. Thus
\[
k r e^{t J} t\left(-J' \right) \leq \frac{k r L z}{\left(1-a\right) \left(1-z\right)} \leq \frac{8 L p^{2}}{1-p^{2}}.
\]
Using $\frac{r}{k}\leq \frac{3}{2}$ now gives
\[
r^{2} q^{2} h' \geq q-\frac{8L p^{2} q}{1+p}-\frac{3p}{2} \log \left(\frac{1}{2p}\right) > \frac{5+15p-102p^{2} +112p^{3}}{20\left(1+p\right)}.
\]
The last inequality uses $\log (\frac{1}{2p}) \leq \frac{1-2p}{2p}$ and $L<\frac{7}{10}$. With
$s=1-2p$ the numerator is $1+3\frac{s}{2}+33\frac{s^{2}}{2}-14s^{3} \geq
1+3\frac{s}{2}+5\frac{s^{2}}{2}>0$.

For the claim about $g$, write $\delta =-\log p-L \theta$. At fixed $k$, put $x_{0} =\log
(\frac{a}{z})$ and $L_{0} (z) =-\log (1-z)$. Then
\[
\begin{aligned}
 \theta&=\frac{a(1-z)L_0(z)}{J(1-a)z},\\
 \frac{\mathrm d\theta}{\mathrm dx_0}
 &=\frac{a}{J(1-a)}\left(\frac{L_0(z)}{z}-1\right)
 \leq B_k:=\frac{J-a}{(1-a)J}.
\end{aligned}
\]
At $x_{0} =0$ we have $\theta =1$. The logarithmic series gives $J-a\leq
\frac{a^{2}}{2}+\frac{a^{3}}{3(1-a)}$ and $J\geq a+\frac{a^{2}}{2}+\frac{a^{3}}{3}$. Hence
\[
B_{k} \leq \frac{a\left(3-a\right)}{6\left(1-a\right)^{2} \left(1+\frac{a}{2}+\frac{a^{2}}{3}\right)}\leq 32\frac{a}{45}.
\]
The rational factor multiplying $a$ is increasing for $a\leq \frac{1}{4}$, as its derivative has the
sign of $21-12a+11a^{2} -26a^{3} +6a^{4} >0$. Since $k a\leq \frac{1}{2}$, $k L B_{k}
<\frac{56}{225}<\frac{1}{4}$. Integrating from $x_{0} =0$ gives $\theta -1\leq k B_{k} \log
(\frac{1}{2p})$. Consequently $\delta \geq (1-k L B_{k} ) \log (\frac{1}{2p}) >3\frac{1-2p}{4}$ for
$p<\frac{1}{2}$.

Differentiating $g$ shows that $g' \leq 0$ is equivalent to
\[
\frac{r}{k} \left[k-u^{2} \left(1-z\right) \right] \delta \geq \left(1-2p\right) \left[1+u\left(1-z\right) -\frac{k r L \theta z}{1-z}\right] .
\]
In view of the strict bound for $\delta$, $\theta \geq 1$ and $u\leq 1$, it suffices to show
\[
R\left(k,p\right) :=1+u\left(1-p^{k} \right) -\frac{k r L p^{k}}{1-p^{k}}\leq M\left(k\right) :=\frac{3}{4}\left(k-\frac{1}{k}\right) .
\]
At $k=2$, $R\leq 1+\frac{p}{q}-4p^{2} =1+\frac{p(1-2p)^{2}}{q}\leq 1+p(1-2p) \leq \frac{9}{8}=M(2)$.
For the derivative at fixed $p$ put $b_{0} =-\log p$ and $x=k b_{0}$. With $A(x) =(x-2) e^{x} +2$
and $B(x) =(x-1) e^{x} +1$, direct differentiation gives
\[
\partial_{k} R=u p^{k} b_{0} +L \frac{k A\left(x\right) +B\left(x\right)}{\left(e^{x} -1\right)^{2}}.
\]
The first term is at most $\frac{L}{4}<\frac{7}{40}$, using $\frac{-p \log p}{q}\leq L$. For the
second, $B(x) \geq 0$ and $2B(x) \leq (e^{x} -1)^{2}$ by $\sinh x\geq x$. As $k\leq \frac{x}{L}$,
the term is bounded by $\frac{1}{2}$ provided
\[
2x A\left(x\right) +\frac{7}{5} B\left(x\right) \leq \left(e^{x} -1\right)^{2} .
\]
The difference has zero constant and linear coefficients. Its coefficient at $x^{n}$ for $n\geq 2$
is $\frac{5 \cdot 2^{n} -10n^{2} +23n-3}{5n!}$. The numerators at two and three are 23 and 16. For
$n=m+4$ the bound $2^{m} \geq 1+m+\frac{m(m-1)}{2}$ leaves at least $30m^{2} -17m+9>0$. Thus
$\partial_{k} R<\frac{27}{40}<\frac{3}{4}\leq M' (k)$. This proves the comparison and the strict
decrease of $g$ for $0<t<1$. \end{proof}

We now compare the integrands at a fixed quantile. Under the hypotheses on $H$ below,
the three monotonicity statements just proved imply that the denominator is
nonincreasing in $k$.

\begin{lemma}\label{r4-6}
Use $\mathcal{K}_{d} =d \alpha_{k} +\omega_{k} H(p)$ in the integral defining $c_{d}^{\ast}$.
Suppose $H$ is positive and convex on $[0,\frac{1}{2}]$, $\frac{H(p)}{p}$ is nonincreasing and
$H(\frac{1}{2}) \leq \frac{d}{2}$. Then $\frac{c_{d}^{\ast}}{J_{k}}$ is nondecreasing in integer
$k\geq 2$. These conditions hold for the functions $\widehat{M}_{d}$ in \cref{r4-2} when $d\geq 3$.
\end{lemma}
\begin{proof} Put $F_{k} =\frac{1-p^{k}}{1-p}$ and $w_{k} =k \frac{p^{k-1}}{1-p^{k}}$. The
definitions give $\alpha_{k} +k p^{k} =p(r F_{k} -\omega_{k} )$ and $h_{k} =\frac{\omega_{k}}{r
F_{k}}$. With $s(p) =1-\frac{H(p)}{d p}$ we get
\[
\frac{d c_{d}^{\ast}}{J_{k}}=\frac{1}{J_{k}} \int_{0}^{1/2} \frac{w_{k}}{1-h_{k} s} \,\mathrm{d} p=\int_{0}^{1} \frac{\mathrm{d} t}{1-h_{k} \left(p_{k} \left(t\right) \right) s\left(p_{k} \left(t\right) \right)}.
\]
The denominators are positive since $0<h_{k} <1$ and $s<1$. Write $s_{\ast} =s(\frac{1}{2}) \geq 0$ and
\[
C\left(p\right) =\frac{H\left(p\right) -2p H\left(\frac{1}{2}\right)}{d\left(1-2p\right)}=\frac{1}{d} \left[H\left(\frac{1}{2}\right) -\frac{1}{2} \frac{H\left(\frac{1}{2}\right) -H\left(p\right)}{\frac{1}{2}-p}\right] .
\]
Then $C\geq 0$ by monotonicity of $\frac{H(p)}{p}$, and $C$ is nonincreasing by convexity of $H$.
Since $h_{k} s=s_{\ast} h_{k} -g C$, \cref{r4-5} shows that $h_{k} s$ is nondecreasing at each
quantile. Integrating proves the result. This identity for the transformed integral and the
decomposition of $h_{k} s$ also hold for positive $H$ with $s_{\ast} <0$.

For  $\widehat{M}_{3} ,\widehat{M}_{4} ,M_{d-1}$ the ratio $\frac{H}{p}$ is nonincreasing by
differentiation, and $H(\frac{1}{2}) \leq \frac{d}{2}$ (for $M_{d-1}$ this reduces to $d\leq
2^{d-1}$). Convexity is immediate for $\widehat{M}_{3} ,\widehat{M}_{4}$. For $m=d-1$, the second
derivative of $M_{m}$ equals
\[
\frac{m^{2} q^{m-2}}{\left(1-q^{m} \right)^{3}} \left[\left(m-1\right) -\left(m+1\right) q+\left(m+1\right) q^{m} -\left(m-1\right) q^{m+1} \right] .
\]
The bracket vanishes at $q=1$. Its derivative is
\[
\left(m+1\right) \left(-1+m q^{m-1} -\left(m-1\right) q^{m} \right) \leq 0,
\]
since $m q^{m-1} -(m-1) q^{m}$ increases to one. This proves convexity, including zero by continuity. \end{proof}

\subsection{Large degrees}\label{large}
For $d\geq 10$ the weights are linear, since $\kappa =0$. By \cref{r4-4},
the first assertion of Proposition~\ref{degree} follows once $c_{d}^{\ast} >c$.

\begin{lemma}\label{r4-7}
For $d\geq 10$ and $k\geq 2$, $c_{d}^{\ast} >\frac{J_{k} (1+\frac{1}{10d})}{d}$.
\end{lemma}
\begin{proof} By \cref{r4-6} it suffices to take $k=2$. Thus $J=J_{\ast}$, $\mathcal{K}_{d} =d
p+(2+p) M_{d-1}$ and $c_{d}^{\ast} =6 \int_{0}^{1/2} \frac{p^{2}}{q(\mathcal{K}_{d} +2d p^{2}
)}\,\mathrm{d} p$. Write
\[
w=\frac{2p}{1-p^{2}}, \quad h=\frac{2+p}{3\left(1+p\right)}, \quad Q_{0} =q^{d-1} , \quad u=\frac{h}{d}+\left(1-h\right) Q_{0} .
\]
All integrals in this proof are over $[0,\frac{1}{2}]$ unless specified. Substitution gives $d
c_{d}^{\ast} =\int \frac{w(1-Q_{0} )}{1-u} \,\mathrm{d} p$ and $\int w \,\mathrm{d} p=J$. Set $A=w
h$, $B=w(1-h)$, $T=\int w Q_{0} \,\mathrm{d} p$ and
\[
L=\int A \,\mathrm{d} p-\left(d+1\right) \int A Q_{0} \,\mathrm{d} p-d \int B Q_{0}^{2} \,\mathrm{d} p.
\]
Cauchy's inequality gives $d c_{d}^{\ast} \geq \frac{(J-T)^{2}}{J-2T-\frac{L}{d}}$, with positive
denominator $\int w(1-Q_{0} ) (1-u) \,\mathrm{d} p$. On clearing it, the desired strict inequality
follows from $L>0$ and $L+d \frac{T^{2}}{J}+\frac{T}{5}>\frac{J}{10}$.

Put $P_{A} =4\frac{p}{3}-2\frac{p^{2}}{3}+2p^{3}$, $P_{B}
=2\frac{p}{3}+2\frac{p^{2}}{3}+3\frac{p^{4}}{2}$, and $D_{0} =(1+p)^{2} (1-p)$. The exact
differences are
\[
\begin{aligned}
 D_0(P_A-A)&=\tfrac{2}{3}p^4(2-2p-3p^2),\\
 D_0(P_B-B)&=p^4\left(\tfrac16+\tfrac56 p-\tfrac32 p^2-\tfrac32 p^3\right).
\end{aligned}
\]
Both are nonnegative on $[0,\frac{1}{2}]$ by monotonicity of the first bracket and concavity of the
second. Since $D_{0} \geq 1$, integrating $P_{A} -(\frac{2}{3}) p^{4} (2-2p-3p^{2} )$ gives $\int A
\,\mathrm{d} p\geq \frac{3377}{20160}>\frac{1}{6}$. Both polynomials are nonnegative on $[0,1]$.
Extend their weighted integrals to that interval and integrate their monomials to obtain
\[
\left(d+1\right) \int A Q_{0} \,\mathrm{d} p\leq E_{A} :=\frac{4\left(d^{2} +4d+12\right)}{3d\left(d+2\right) \left(d+3\right)},
\]
\[
\begin{aligned}
 d\int BQ_0^2\,\mathrm dp &\leq E_B :=
 \frac{1}{3(2d-1)}+\frac{2}{3(2d-1)(2d+1)}\\
 &\hphantom{\leq E_B :={}}+\frac{18}{(2d-1)(2d+1)(2d+2)(2d+3)}.
\end{aligned}
\]
Both bounds decrease in $d$. In particular $-E_{A}' =4\frac{d^{4} +8d^{3} +50d^{2} +120d+72}{3d^{2}
(d+2)^{2} (d+3)^{2}}>\frac{4}{3(d+3)^{2}}$. Since $E_{A} (10) <\frac{13}{100}$ and $E_{B} (10)
<\frac{29}{1500}$, we have $L>\frac{1}{6}-E_{A} -E_{B} >\frac{13}{750}>0$. Moreover $w>2p$ on the
interior, whence
\[
T> \frac{2\left(1-\left(d+2\right) 2^{-d-1} \right)}{d\left(d+1\right)}>\frac{99}{50d\left(d+1\right)}.
\]
The factor $(d+2) 2^{-d-1}$ decreases from $\frac{3}{512}$ at ten. The series for $2
\operatorname{arctanh}(\frac{1}{7})$ gives $J<\frac{145}{504}<\frac{36}{125}$. The estimates yield
$L+d \frac{T^{2}}{J}+\frac{T}{5}>R(d)$, where
\[
R\left(d\right) =\frac{1}{6}-E_{A} -E_{B} +\frac{1089}{80d\left(d+1\right)^{2}}+\frac{99}{250d\left(d+1\right)}.
\]
Since $d+3\leq 13\frac{d}{10}$, direct differentiation gives
\[
R' \left(d\right) >\frac{4}{3\left(d+3\right)^{2}}-\frac{3267}{80d^{4}}-\frac{99}{125d^{3}} > \frac{1}{d^{2}} \left(\frac{3}{4}-\frac{41}{100}-\frac{2}{25}\right) >0.
\]
At $d=10$, $R(10)
>\frac{13}{750}+\frac{9}{800}+\frac{9}{2500}=\frac{1931}{60000}>\frac{18}{625}>\frac{J}{10}$. This
proves the assertion. \end{proof}

Lemma~\ref{r4-7} gives the bound on the remainder for linear weights.
In the reconstruction inequality of Proposition~\ref{degree}, we need strict
inequality at every positive residual degree to rule out further cases of equality.

\begin{lemma}\label{r4-8}
For $d\geq 10$, the linear weights satisfy the second part of \cref{degree}.
\end{lemma}
\begin{proof} Set $v=2^{-d}$, $\epsilon =\frac{J}{10d}$ and $a=J+\epsilon =d c$. For $0\leq t\leq d$, put
\[
C\left(t\right) =1+v e^{\epsilon t} \left(e^{d J} -e^{J t} \right) , \quad \delta =e^{\frac{J}{10}} v\left(e^{d J} -1\right) .
\]
Then $Q^{d-t} +(1-Q^{d-t} ) v e^{a t} =Q^{d} e^{J t} C(t)$ and $b=Q^{d} C(0)$. After taking
logarithms it suffices, by integration from zero, to show $r d(\log C)' (t) +k \log C(0)
<\frac{J}{10}$. Since $1\leq C\leq 1+\delta$ and $C' =\epsilon (C-1) -J v e^{a t}$, we have $(\log
C)' \leq \epsilon \delta -J \frac{v}{1+\delta}$. Put $H(w) =\frac{e^{w} -1}{w}$. Using $\log C(0)
\leq d v J H(d J)$ gives
\[
\frac{r d\left(\log C\right)' + k \log C\left(0\right)}{J} \leq r \delta \left(\frac{1}{10}+d v\right) +d v k\left(H\left(d J\right) -1\right) -d v.
\]
The right side is largest at $k=2$. Indeed $J_{k+1} <\frac{J_{k}}{2}$ by strict convexity of $-\log
(1-x)$. Also $\frac{\delta}{J}=v \int_{0}^{d} e^{J(s+\frac{1}{10})} \,\mathrm{d} s$ and $\frac{H(d
J) -1}{J}=\int_{0}^{1} \frac{e^{d J z} -1}{J} \,\mathrm{d} z$ increase with $J$. Hence $r \delta$
and $k(H(d J) -1)$ decrease with $k$.

At $k=2$, let $j=\log (\frac{4}{3})$, $A=e^{\frac{j}{10}}$ and $u=(\frac{2}{3})^{d}$. The bound
becomes $(\frac{2}{j}+3\frac{A}{10}) (u-v) -3d v+3A d v(u-v)$. With $\frac{2}{7}<j<\frac{1}{3}$ and
$A<\frac{30}{29}$ this is less than
\[
R_{d} =\frac{212}{29} \left(u-v\right) -3d v+\frac{90}{29} d v\left(u-v\right) .
\]
Here $(\frac{4}{3})^{d} \geq d+1$ for $d\geq 10$, by induction from ten. Put $P_{d} =(\frac{212}{29}) (u-v) -3d v$. Then
\[
P_{d+1} -P_{d} =\left(\frac{212}{29}\right) \left(-\frac{u}{3}+\frac{v}{2}\right) +3\left(d-1\right) \frac{v}{2}< -\left(d+1\right) \frac{v}{2}.
\]
Also $W_{d} =d v(u-v)$ decreases, since
\[
W_{d+1} -W_{d} =v^{2} \left[-\left(2d-1\right) \frac{u}{3v}+\frac{3d-1}{4}\right] <0.
\] Thus $R_{d} \leq R_{10}$. Using $\frac{212}{29}<\frac{22}{3}$, $\frac{90}{29}<\frac{25}{8}$ and $(\frac{2}{3})^{10} <\frac{7}{400}$ gives
\[
R_{10} <\left(\frac{22}{3}+\frac{125}{4096}\right) \frac{7}{400}-\frac{15}{512}=\frac{489409}{4915200}<\frac{1}{10}.
\]
This proves the strict comparison for $t>0$ and equality at zero. Finally $\rho_{\ell} =c \ell \leq
J(1+\frac{1}{10d}) <\frac{101}{300}<\log 2$. \end{proof}

\subsection{Small degrees}\label{small}
Suppose $2\leq d\leq 9$. Write $v=2^{-d}$, $c_{0} =-\frac{\log b}{d^{2}}$, $\tau =3\frac{4d+1}{20d}$
and $\sigma =\frac{J^{2}}{9J_{\ast}}$. Then $c=c_{0} +r \frac{\sigma}{10d^{2}}$ and $\kappa =\tau v
\sigma$. Dividing \eqref{target} by $\frac{\ell}{r}$ gives
\begin{equation}\label{smalltarget}
c_{d}^{\ast} -c+\frac{r}{\ell} \Delta_{d,\ell} -r \kappa \left(d-\ell \right) >0 \quad \left(1\leq \ell \leq d\right) .
\end{equation}
For $k=2,3,4$ we check this inequality by rational arithmetic. For larger $k$
we compare the integral with a bound at rank five.

For a polynomial $P=\sum_{j} P_{j} p^{j}$, put
\[
\mathcal{I}\left(P\right) =\sum_{j} \frac{P_{j}}{\left(j+1\right) 2^{j+1}}, \quad \mathcal{C}\left(P,D\right) =\frac{\mathcal{I}\left(P\right)^{2}}{\mathcal{I}\left(P D\right)}.
\]
Thus $\mathcal{I}$ integrates over $[0,\frac{1}{2}]$. Cauchy's inequality gives $\mathcal{C}(P,D)
\leq \int_{0}^{1/2} \frac{P}{D} \,\mathrm{d} p$ for $P\geq 0$, $P$ nonzero and $D>0$ in the interior
(with finite integrals). We also use the rational bounds $L(x) \leq \log x\leq U(x)$ for $x\geq 1$.
With $z=\frac{x-1}{x+1}$ they are
\[
L\left(x\right) =2\left(z+\frac{z^{3}}{3}\right) , \quad U\left(x\right) =L\left(x\right) +\frac{2z^{5}}{5\left(1-z^{2} \right)}.
\]
They follow from the series for $2 \operatorname{arctanh} z$. Set $J_{+} =U(\frac{1}{Q})$, $J_{\ast
,-} =L(\frac{4}{3}) =\frac{296}{1029}$ and $\sigma_{+} =\frac{J_{+}^{2}}{9J_{\ast ,-}}$.

\begin{lemma}\label{r4-9}
The comparison \eqref{smalltarget} holds for $2\leq d\leq 9$ and every $k\geq 2$.
\end{lemma}
\begin{proof} Let $H_{m} =\sum_{j=0}^{m-1} q^{j}$, $A=\alpha_{k} +k p^{k}$, $T=q^{d-1} (1-p^{2}
)^{2}$ and $S=1-p^{k} +p^{k} T$. The identity $\mathcal{K}_{d} +k d p^{k} =\frac{N}{h}$ uses the
following polynomials:
\begin{center}\small\begin{tabular}{lcl}
\toprule
$d$ & $h$ & $N$ \\
\midrule
$2$ & $1-\frac{p^{2}}{2}$ & $2A h+\omega_{k} -b_{k} q^{2}$ \\
$3$ & $3-2p$ & $3A h+3\omega_{k}$ \\
$4$ & $(2-p)^{2}$ & $4A h+4\omega_{k}$ \\
$d\geq 5$ & $H_{d-1}$ & $d A h+(d-1) \omega_{k}$ \\
\bottomrule
\end{tabular}\end{center}
For $\ell <d$ set $V_{\ell} =(N-\ell A h) H_{\ell} -\ell \omega_{k} h$. Dividing by $h H_{\ell}$
gives $\mathcal{K}_{d} -\mathcal{U}_{k,\ell} +k(d-\ell ) p^{k} >0$ in the interior. The defining
integrals and Cauchy's inequality now give
\[
c_{L} =r k \mathcal{C}\left(p^{k} h,q N\right) \leq c_{d}^{\ast} , \quad D_{\ell}^{L} =\mathcal{C}\left(V_{\ell} p^{k} q^{d-2} \left(1-p^{2} \right)^{2} ,N S\right) \leq \Delta_{d,\ell} .
\]
Here $D_{d}^{L} =0$. The integrand for $\Delta_{d,\ell}$ is the first polynomial argument divided by the second.

\emph{The ranks $k=2,3,4$.} Since $b=Q^{d} C_{d}$ for $C_{d} =1-v+v Q^{-d}$, upper bounds for $c,\kappa$ are
\[
c_{+} =\frac{J_{+}}{d}-\frac{L\left(C_{d} \right)}{d^{2}}+\frac{r \sigma_{+}}{10d^{2}}, \quad \kappa_{+} =\tau v \sigma_{+} .
\]
For each positive degree $\ell$, put
\[
 G_\ell=c_L-c_+ +\frac r\ell D_\ell^L-r\kappa_+(d-\ell).
\]
For each pair $(k,d)$ the table gives the integer part of
$10^6\min_{1\leq\ell\leq d}(G_\ell/J_+)$. Its entries follow by rational
substitution in the preceding formulas.
\begin{center}\small\begin{tabular}{rrrrrrrrr}
\toprule
$d$ & 2 & 3 & 4 & 5 & 6 & 7 & 8 & 9 \\
\midrule
$k=2$ & 24584 & 394 & 2412 & 2623 & 2180 & 1964 & 1791 & 1548 \\
$k=3$ & 30603 & 5187 & 10922 & 9517 & 7654 & 6477 & 5585 & 4850 \\
$k=4$ & 36291 & 10332 & 16301 & 15019 & 11763 & 9585 & 7985 & 6738 \\
\bottomrule
\end{tabular}\end{center}
For example, at $(k,d,\ell ) =(2,3,3)$ one obtains $c_{L} =\frac{45}{566}$, $c_{+}
<\frac{397}{5000}$, and $G_{3} >\frac{149}{1415000}$. The table covers all $3 \sum_{d=2}^{9} d=132$
comparisons.

\emph{The ranks $k\geq 5$ with $3\leq d\leq 9$.} In this range the edge rate alone suffices.
Concavity of the logarithm gives $\frac{c_{0}}{J}\leq \frac{1-v}{d}$. Put $\gamma
=\frac{1}{10d^{2}}+\tau v(d-1)$. Then
\[
\frac{c+r \kappa \left(d-1\right)}{J} \leq \frac{1-v}{d}+\frac{r J}{9J_{\ast}} \gamma .
\]
The right side decreases with $k$ because $J_{k+1} <\frac{J_{k}}{2}$, while $\frac{c_{d}^{\ast}}{J}$
increases weakly by \cref{r4-6}. At $k=5$, put $J_{5,+} =U(\frac{32}{31})
=\frac{236254723}{7441398720}$ and use the preceding $c_{L}$ at rank five. Exact lower bounds for
\[
10^{6} \left(\frac{c_{L}}{J_{5,+}}-\frac{1-v}{d}-\frac{6J_{5,+}}{9J_{\ast ,-}} \gamma \right)
\]
are, in order of $d=3,\dots ,9$,
\[
914, \quad 10825, \quad 12351, \quad 10710, \quad 9418, \quad 8256, \quad 7198.
\]
For $d=3$, $c_{L} =\frac{3250}{335237}$. Hence $c_{d}^{\ast} >c+r \kappa (d-1)$ throughout this range, as required.

\emph{The ranks $k\geq 5$ with $d=2$.} Replacing the two-edge score by $2 \alpha_{k} +\omega_{k} H$,
$H=(1-\frac{p^{2}}{2})^{-1}$, increases the score bound and lowers the integrated rate. With the
notation of \cref{r4-6}, $s=1-\frac{H}{2p}$, $s_{\ast} =-\frac{1}{7}$, and $C=\frac{H(p) -2p
H(\frac{1}{2})}{2(1-2p)}$ is nonnegative and nonincreasing. At the same quantile and for $k\geq K$,
\cref{r4-5} gives $g_{k} \leq g_{K}$, $C(p_{k} ) \leq C(p_{K} )$, and $h_{k} \leq 1$. Thus $h_{k}
s(p_{k} ) \geq s_{\ast} -g_{K} C(p_{K} )$. Returning to the integral at rank $K$ yields
\[
\frac{c_{2}^{\ast}}{J_{k}} \geq \frac{K\left(K+1\right)}{J_{K}} \int_{0}^{1/2} \frac{p^{K}}{q\left[\left(\frac{16}{7}\right) \left(\alpha_{K} +K p^{K} \right) +\omega_{K} H\right]} \,\mathrm{d} p.
\]
Indeed the denominator in quantile coordinates becomes
\[
1-s_{\ast} +h_{K} \left(s_{\ast} -s\left(p\right) \right) =\frac{8}{7}+h_{K} \left(\frac{H}{2p}-\frac{8}{7}\right).
\] At $K=5$ put
\[
\begin{aligned}
 P&=p^5(1-p^2/2),\\
 D&=q\left[\frac{16}{7}\left(\sum_{s=1}^5 s p^s\right)(1-p^2/2)
              +5+4p+3p^2+2p^3+p^4\right].
\end{aligned}
\]
Now $\mathcal{I}(P) =\frac{29}{12288}$, $\mathcal{I}(P D) =\frac{1115124719}{82657935360}$, and
$30\mathcal{C}(P,D) =\frac{441966525}{35683991008}>\frac{99}{8000}$. The series bound $-\log (1-x)
\leq x+\frac{x^{2}}{2(1-x)}$ gives $J_{5} \leq \frac{63}{1984}$. Thus
$\frac{c_{2}^{\ast}}{J_{k}}>\frac{341}{875}$. Concavity also gives $-\log b\leq 3\frac{J_{k}}{2}$,
so
\[
\frac{c+r \kappa}{J_{k}} \leq \frac{3}{8}+\frac{r J_{k}}{9J_{\ast}} \frac{31}{160} \leq \frac{3987}{10240},
\]
using the decrease of $r J_{k}$ and $J_{\ast} >\frac{2}{7}$. The difference
$\frac{341}{875}-\frac{3987}{10240}=\frac{643}{1792000}$ is positive. This completes all cases.
\end{proof}

It remains to verify the reconstruction inequality for the quadratic weights.
For each fixed residual degree $\ell$ with $1\leq\ell<d$, we reduce the required
estimate to $k=2$.

\begin{lemma}\label{r4-10}
For $2\leq d\leq 9$ the quadratic weights satisfy the second part of \cref{degree}.
\end{lemma}
\begin{proof} Write $j=-\frac{\log b}{d}$ and $P_{\ell} =\frac{\ell}{10d^{2}}+\tau v \ell (d-\ell )$. Then
\[
\phi_{\ell} =\frac{j \ell}{r d}+\sigma P_{\ell} , \quad \rho_{\ell} =\frac{j \ell}{d}+\sigma \left(P_{\ell} +k \ell P_{1} \right) .
\]
We have $P_{\ell} \leq \ell P_{1}$ and $P_{1} \leq \frac{1}{40}+\frac{27}{160}<\frac{1}{5}$, since
$\tau \leq \frac{27}{40}$ and $(d-1) 2^{-d} \leq \frac{1}{4}$. Since $J_{k+1} <\frac{J_{k}}{2}$,
$(k+1)^{2} J_{k}$ decreases for $k\geq 2$ and $\sigma \leq \frac{J}{r^{2}}$. As $0\leq j\leq J$, we
obtain $\rho_{\ell} \leq J+J \ell \frac{P_{1}}{r}<8\frac{J}{5}<\log 2$ for $1\leq \ell \leq d$. Also
$\rho_{0} =0$.

At $\ell =0$ the reconstruction comparison is equality. At $\ell =d$ the left side is zero and the
right side is $d \sigma P_{d} >0$. For $1\leq \ell <d$ put $s=d-\ell$, $a=\frac{\ell}{d}$, and
\[
\begin{aligned}
 F(z)&=\log(1-v+ve^z),\qquad D(J)=aF(dJ)-F(\ell J),\\
 \beta(J)&=v(e^{dJ}-e^{\ell J})e^{-aF(dJ)} .
\end{aligned}
\]
With $u=d \sigma P_{\ell}$, $w=d k \sigma \ell P_{1}$, $A=e^{-s(J-j)}$ and $B=(\frac{v}{b}) (1-e^{-s
J} )$, the desired inequality is $A e^{-u} +B e^{w} \leq 1$. Since $J-j=\frac{F(d J)}{d}$, we have
$\frac{B}{A}=\beta (J)$ and $\frac{A+B-1}{A}=1-e^{-D(J)}$. Convexity gives $D(J) \geq 0$. It
therefore suffices to prove
\begin{equation}\label{reconbase}
D\left(J\right) +\beta \left(J\right) \left(e^{w} -1\right) <1-e^{-u} .
\end{equation}
We reduce \eqref{reconbase} to $k=2$. For $0\leq z\leq d J_{\ast}$, $\pi (z) =v \frac{e^{z}}{1-v+v
e^{z}}<\frac{1}{2}$, since $v e^{z} \leq (\frac{2}{3})^{d} \leq \frac{4}{9}<1-v$. Thus $F'' =\pi
(1-\pi )$ is increasing. The identity
\[
\frac{a F\left(x\right) -F\left(a x\right)}{x^{2}}=\left(1-a\right) \int_{0}^{a} t F'' \left(t x\right) \,\mathrm{d} t+a \int_{a}^{1} \left(1-t\right) F'' \left(t x\right) \,\mathrm{d} t
\]
shows that $\frac{D(J)}{J^{2}}$ increases with $J$. Likewise $\frac{\beta (J)}{J}=v[\frac{e^{s J}
-1}{J}] e^{\ell J-a F(d J)}$ increases with $J$, since the exponent has derivative $\ell (1-F' (d J)
) >0$.

Put $u_{\ast} =d J_{\ast} \frac{P_{\ell}}{9}$ and $w_{\ast} =2d J_{\ast} \ell \frac{P_{1}}{9}$, the
values at $k=2$. The sequences $J_{k}^{2} ,k J_{k}^{2} ,k J_{k}$ decrease, hence $u\leq u_{\ast}$,
$w\leq w_{\ast}$, $\frac{D(J)}{u}\leq \frac{D(J_{\ast} )}{u_{\ast}}$, and $k \beta (J) \leq 2 \beta
(J_{\ast} )$. Since $\frac{w}{u}=k \ell \frac{P_{1}}{P_{\ell}}$ and $\frac{e^{x} -1}{x}$ increases,
the left side of the sufficient inequality divided by $u$ is at most its value at $k=2$. The right
side divided by $u$ is at least its value there because $\frac{1-e^{-x}}{x}$ decreases. It remains
to verify the strict inequality at $J_{\ast}$.

For $x\geq 0$ put $S(x) =1+x+\frac{x^{2}}{2}$ and, for $x<4$, $T(x) =S(x)
+\frac{x^{3}}{6(1-\frac{x}{4})}$. The exponential series gives $S(x) \leq e^{x} \leq T(x)$ on this
interval. Write $J_{-} =L(\frac{4}{3})$, $J_{+} =U(\frac{4}{3})$ and $C_{m}
=1-v+v(\frac{4}{3})^{m}$. Define
\[
u_{-} =d P_{\ell} \frac{J_{-}}{9}, \quad u_{+} =d P_{\ell} \frac{J_{+}}{9}, \quad w_{+} =2d \ell P_{1} \frac{J_{+}}{9},
\]
\[
D_{+} =a U\left(C_{d} \right) -L\left(C_{\ell} \right) , \quad \beta_{+} =\frac{v\left(\left(\frac{4}{3}\right)^{d} -\left(\frac{4}{3}\right)^{\ell} \right)}{S\left(a L\left(C_{d} \right) \right)}.
\]
Here $J_{+} \leq \frac{7}{24}<\frac{1}{3}$, $P_{1} <\frac{1}{5}$ and $\ell \leq 8$ give $w_{+}
<\frac{16}{15}$. Since $u_{\ast} \geq u_{-} \geq 0$, $w_{\ast} \leq w_{+}$, $D(J_{\ast} ) \leq
D_{+}$ and $\beta (J_{\ast} ) \leq \beta_{+}$, a rational lower bound for the required difference is
\[
\Xi_{d,\ell} =1-\frac{1}{S\left(u_{-} \right)}-D_{+} -\beta_{+} \left(T\left(w_{+} \right) -1\right) .
\]
The table gives $\lfloor 10^{6} \min_{1\leq \ell <d} \frac{\Xi_{d,\ell}}{u_{+}}\rfloor$ by direct rational substitution.
\begin{center}\small\begin{tabular}{rrrrrrrrr}
\toprule
$d$ & 2 & 3 & 4 & 5 & 6 & 7 & 8 & 9 \\
\midrule
Bound & 77536 & 19102 & 8774 & 25859 & 40407 & 55991 & 80799 & 107191 \\
\bottomrule
\end{tabular}\end{center}
All 36 comparisons are strict. This proves the reconstruction bound for every rank. \end{proof}

By \cref{r4-7} and \cref{r4-9}, \eqref{target} holds at every positive degree. \cref{r4-4} gives the
hereditary bound, and \cref{r4-8,r4-10} give the reconstruction comparison and its strictness. This
proves \cref{degree}.
\section*{Acknowledgements}
Caihong Yang was supported by the Scientific Research
Funds at China University of Geosciences (Wuhan) (Project No.2026039). Guorong Gao was supported by the National Key R\&D Program of China (Grant No.~2023YFA1010202), the National Natural Science Foundation of China (Grant No.~12401448), and the Natural Science Foundation of Fujian Province (Grant No.~2024J08030).

\section*{Declaration on the use of generative AI}
The authors used GPT-6 to discuss the choice of local statistics in the proof of Lemma~\ref{r4-4} and the change of variables underlying the comparison across ranks in Lemmas~\ref{r4-5} and~\ref{r4-6}. It also helped us identify and apply the weighted bipartite graph inequality in Lemma~\ref{r3-2}, check parts of the proofs and the calculations for small degrees. 
\bibliographystyle{abbrv}
\bibliography{independent-sets}

\end{document}